\documentclass[a4paper,10pt]{amsart}
\usepackage[ansinew]{inputenc}
 
\usepackage{geometry}

\usepackage{amsmath}
\usepackage{amsfonts}
\usepackage{amsthm}
\usepackage{amssymb}
\usepackage{graphicx}
\makeindex

\usepackage{array} %centerd table cells m{}m{}
\usepackage[T1]{fontenc} %small caps
\newcommand{\subtile}[1] %subtitle
{
	\vspace{-0.3cm}
	\begin{center}
 		{{\textsc{#1}}}\\
	\end{center}
	\vspace{0.1cm}
}

\newcommand{\imageplain}[2] %Image scaled
{
      {\includegraphics[scale=#1]{#2.png}}
}
\newcommand{\image}[2] %Image scaled centered
{
      \begin{center} %Image
	\includegraphics[scale=#1]{#2.png}
      \end{center}
}
\newcommand{\hamburger}[4] 
{
  \thispagestyle{empty}
  \vspace*{-2cm}
  \begin{flushright}
    ZMP-HH / #2 \\
    Hamburger Beitr{\"a}ge zur Mathematik Nr. #3 \\
   #4 \\
  \end{flushright}
  \vspace{0.5cm}
  \begin{center}
    \Large \bf
    #1
  \end{center}
  \vspace{0.5cm}
  \begin{center}	
    Simon Lentner\\
    Algebra and Number Theory, 
    University Hamburg,\\
    Bundesstra{\ss}e 55, D-20146 Hamburg \\
  \end{center}
  \vspace{0.5cm}

}

\makeatletter
\newcommand{\xleftrightarrow}[2][]{\ext@arrow 3359\leftrightarrowfill@{#1}{#2}}
\newcommand{\xdashrightarrow}[2][]{\ext@arrow 0359\rightarrowfill@@{#1}{#2}}
\newcommand{\xdashleftarrow}[2][]{\ext@arrow 3095\leftarrowfill@@{#1}{#2}}
\newcommand{\xdashleftrightarrow}[2][]{\ext@arrow
3359\leftrightarrowfill@@{#1}{#2}}
\def\rightarrowfill@@{\arrowfill@@\relax\relbar\rightarrow}
\def\leftarrowfill@@{\arrowfill@@\leftarrow\relbar\relax}
\def\leftrightarrowfill@@{\arrowfill@@\leftarrow\relbar\rightarrow}
\def\arrowfill@@#1#2#3#4{%
  $\m@th\thickmuskip0mu\medmuskip\thickmuskip\thinmuskip\thickmuskip
   \relax#4#1
   \xleaders\hbox{$#4#2$}\hfill
   #3$%
}
\makeatother

\newcommand{\Ksymb}[3]{      
      \begin{bmatrix}
      #1; #2 \\
      #3   
      \end{bmatrix}
}

\allowdisplaybreaks[1]

\theoremstyle{plain}
\newtheorem{theorem}{Theorem}[section]

\newtheorem*{acknowledgementX}{Acknowledgement}

\newtheorem{corollary}[theorem]{Corollary}
\newtheorem{observation}[theorem]{Observation}

\newtheorem{definition}[theorem]{Definition}

\newtheorem{example}[theorem]{Example}

\newtheorem{question}[theorem]{Question}

\newcommand{\g}{\mathfrak{g}}
\newcommand{\B}{\mathfrak{B}}
\newcommand{\ord}{\mathrm{ord}}

\newcommand{\Z}{\mathbb{Z}}
\newcommand{\Q}{\mathbb{Q}}
\newcommand{\C}{\mathbb{C}}
\newcommand{\N}{\mathbb{N}}
\renewcommand{\L}{\mathcal{L}}
\newcommand{\K}{\mathcal{K}}
\newcommand{\h}{\mathfrak{h}}
\renewcommand{\sl}{\mathfrak{sl}}

\newcommand{\res}{\mathrm{Res}}
\newcommand{\Y}{\mathrm{Y}}
\renewcommand{\exp}[1]{\mathrm{e}^{#1}}

\usepackage{enumerate}

\begin{document}

\hamburger{The unrolled quantum group inside\\ Lusztig's quantum group of divided powers}{17-7}{648}{February 2017}
\begin{abstract}
  In this letter we prove that the unrolled small quantum group, appearing in quantum topology, is a Hopf subalgebra of Lusztig's quantum group of divided powers. We do so by writing down non-obvious primitive elements with the correct adjoint action.
  
  As application we explain how this gives a realization of the unrolled quantum group as operators on a conformal field theory and match some calculations on this side. In particular our results explain a prominent weight shift that appears in \cite{FT10}.

  Our result extends to other Nichols algebras of diagonal type, including super Lie algebras.   
\end{abstract}
\title{}
\date{}
\maketitle
\tableofcontents

\section{Introduction}

Unrolled quantum groups are certain Hopf algebras that were introduced in the context of quantum topology \cite{GPT09,CGP15,GP16}. They are used to construct topological invariants from non-semisimple tensor categories, more precisely the representation category of versions of quantum groups at a root of unity $q$ of even order $\ell$. This approach is in particular able to recover the Reshetikhin-Turaev invariant of $3$-manifolds \cite{RT91}. \\

The idea of ``unrolling'' addresses the following problem: For modules over the small quantum group $u_q(\g)$  at an $\ell$-th root of unity, the action of the group elements $K_\alpha$ in the Cartan part cannot distinuish weight spaces $V_\lambda,V_\mu$, for which $\lambda-\mu$ is a multiple of $\ell$ of fundamental coweights. For even $\ell$ this effect causes additional problems: The braiding factor $q^{(\mu,\lambda)}$ is not a well-defined bimultiplicative map on the finite group of weights modulo $\ell$. Indeed, for even $\ell$ the small quantum group does not necessarily admit a braiding \cite{KS11}.

On way to remedy this is to introduce an \emph{unrolled quantum group}, which is roughly the semidirect product of the quantum group with the Cartan part of the corresponding Lie algebra, which acts differently on each weight space (an alternative remedy is to accept the appearance of nontrivial associators, see Section 4). This theory of unrolling has been successfully worked out for $\sl_2$ in its quantum topology application, it has been generalized to higher rank Lie algebras in \cite{GP16} and to Nichols algebras of diagonal type in \cite{AS17}. \\ 

The main result of this letter is to prove that the unrolled small quantum group appears as a Hopf subalgebra in Lusztig's quantum group of divided powers $U_q^\L(\g)$. This is an infinite-dimensional Hopf algebras constructed by specialization, which admits a Hopf algebra homomorphism
$$U_q^\L(\g)\stackrel{\mathrm{Frob}}{\longrightarrow} U(\g^{(\ell)})$$ 
In the commonly known case  $2\nmid\ell$ (and $3\nmid \ell$ for $G_2$) this result holds for $\g^{(\ell)}=\g$ \cite{Lus90,A96}, while in divisible cases we have $\g^{(\ell)}=\g^\vee$ the Langlands dual Lie algebra, or a variant thereof including an additional symmetric braiding, as discussed in the end of Section 2. \\

From a mathematical perspective, realizing the unrolled small quantum group $u_q(\g)\rtimes U(\h)$ inside $U_q^\L(\g)$ means to construct for this sequence a good section of the Cartan part $U(\g^{(\ell)})^0=U(\h)$. This means to find (non-obvious) preimages $H_\alpha\in U_q^\L(\g)$ of the Chevalley generators of $\mathfrak{h}$, such that the $H_\alpha$ are primitive elements and act on $u_q(\g)$ in the designated way.\\ 

From a physical perspective, part of the Lusztig quantum group $U_q^\L(\g)$ is expected to act by screening charges on a free field theory, more precisely on the lattice vertex algebra associated to a rescaled root lattice $\Lambda$ of $\g$, see e.g. \cite{FGSTsl2,NT11,Len17}. Realizing the unrolled small quantum group inside $U_q^\L(\g)$ means to realize it as operators on this conformal field theory. 

We can indeed present $H_\alpha$ as $\Lambda^*$-grading operators, whereas the group elements $K_\alpha\in u_q(\g)$ act as exponentiated grading operators. The grading operators $H_\alpha$ literally unroll the vertex algebra module corresponding to a $\Lambda$-coset in $\Lambda^*$ into all its $\Lambda^*$-graded subspace, while $K_\alpha$ acts on the entire module as a single scalar. It is widely expected that a vertex subalgebra of this lattice vertex algebra (kernel of short screenings) should have a non-semisimple representation theory equivalent to $u_q(\g)$-representations, so these two unrolling constructions should be closely related.\\

% Besides this main result we construct a curious Hopf algebra $U_q^{\K\L}(\g)$ that acts as a ``hybrid'': It contains the Kac-DeConcini-Procesi quantum group $U_q^{\K}(\g)$ and surjects to the Lusztig quantum group $U_q^{\L}(\g)$. We can realize the full unrolled quantum group $U_q^{\K}(\g)\rtimes U(\h)$ inside this new Hopf algebra. The representation category of $U_q^{\K}(\g)$ fibres over the points of the Lie group $G$ (resp. $G^\vee$) and becomes non-semisimple over points of the subvariety of non-regular conjugacy classes of $G$. The fibre over the unit is the representation category of the small quantum group. Constructing topological invariants depending on the point by unrolling the entire $U_q^{\K}(\g)$ is very interesting.         \\

Finally we mention how our results should extend to quantum doubles of arbitrary Nichols algebras of diagonal type. As most prominent examples this allows us to cover also quantum groups associated to super-Lie algebras, where divided powers appear only for bosonic root vectors.

\section{Preliminaries}

Let $\g$ be a complex semisimple finite-dimensional Lie algebra and let $\mathfrak{h}$ be its Cartan subalgebra. Fix simple roots $\alpha_1,\ldots, \alpha_n\in\mathfrak{h}^*$, where $n=\mathrm{dim}(\mathfrak{h})$ is the rank of $\g$. We denote the corresponding set of positive roots by $\Phi^+$, the Cartan matrix by $a_{ij}$ and the symmetrized Cartan matrix by  $d_{\alpha_i}a_{ij}=(\alpha_i,\alpha_j)$, where $d_\alpha=(\alpha,\alpha)/2\in\{1,2,3\}$ for all $\alpha\in\Phi$. We will usually not denote the dependence on $\g$. \\

Let $q$ be a primitive $\ell$-th root of unity and denote as usual $q_{\alpha}:=q^{d_\alpha}$ and $\ell_\alpha:=\ord(q_{\alpha}^{2})=\ell/(\ell,2d_\alpha)$. Let $v$ be an indeterminate and denote again $v_{\alpha}:=v^{d_\alpha}$. Already Gau{\ss} has introduced and studied $q$-numbers and $q$-factorials in \cite{G11}, and we will use the following common variant, which is obtained by a minor substitution and is symmetric with respect to $q\leftrightarrow q^{-1}$:
$$
[n]_{v_\alpha}:=\frac{v_\alpha^n-v_\alpha^{-n}}{v_\alpha-v_\alpha^{-1}}
=v_\alpha^{-n+1}\frac{(v_\alpha^2)^n-1}{v_\alpha^2-1}
\qquad
[n]_{v_\alpha}!:= [n]_{v_\alpha}\cdots [1]_{v_\alpha}
\qquad
\begin{bmatrix}n\\k\end{bmatrix}_{v_\alpha}:=\frac{[n]_{v_\alpha}!}{[k]_{v_\alpha}![n-k]_{v_\alpha}!}
$$
After polynomial division all these expressions are Laurent polynomials in $\Z[v,v^{-1}]$ and can be specialized to any nonzero complex number. The crucial property for our purposes is $[\ell_\alpha]_{q_\alpha}=0$. For later use, we also note that $q_\alpha^{\ell_\alpha}\in\{+1,-1\}$ and that for $n\in\N_0$ the following formula holds:
$$[\pm n]_{q_\alpha^{\ell_\alpha}}=\pm(q_\alpha^{\ell_\alpha})^{-n+1}\sum_{s=0}^{n-1}(q_\alpha^{2\ell_\alpha})^s =\pm n(q_\alpha^{\ell_\alpha})^{n-1}$$

Let us mention that for our application we are usually interested in the fully divisible case $2d_\alpha\;|\;\ell$, while in quantum group literature (in particular \cite{Lus90} from Sec. 8.4 on) it is usually required that $(2d_\alpha,\ell)=1$. In the general case the orders $\ell_\alpha$ may be different from one another and different from $\ell$. This can cause unusual behaviour, for example the Langlands dual root system $\g^\vee$ appears in the exact sequence below and there are problems with the existence of the standard $R$-matrix. The algebra structure for general $\ell$ are discussed in \cite{Lus93,Len14} and alternative $R$-matrices are determined in \cite{LO16,GLO18}.\\

To avoid degeneracies in this article we exclude small values $\ell\neq 1,2$, and also $\ell\neq 4$ if some $d_\alpha=2$, and $\ell\neq 3,4,6$ if some $d_\alpha=3$; the degenerated cases are treated in \cite{Len14}.\\

To this data the following Hopf algebras are associated by \cite{Lus90}: 
\begin{itemize}
 \item \textbf{Rational form} $U_v^{\Q(v)}(\g)$: An infinite-dimensional Hopf algebra  over the field of rational functions in an indeterminate $\Q(v)$ defined by Drinfeld and Jimbo. As an algebra, it is generated by elements $K_\lambda$ indexed by elements in the root lattice $\lambda\in\Lambda$ (or another lattice between root- and weight-lattice) with $K_{\lambda+\mu}=K_\lambda K_\mu=K_\mu K_\lambda$, together with elements $E_{\alpha_i},F_{\alpha_i}$ for each simple root, and the following relations:
 \begin{align*}
    K_{\lambda}E_{\alpha_i}K_{\lambda}^{-1}
    &=v^{(\lambda,\alpha_i)}E_{\alpha_i}
    ,\; \forall\lambda\in\Lambda
    \qquad \mbox{\emph{(group action)}}\\
    K_{\lambda}F_{\alpha_i}K_{\lambda}^{-1}
    &=v^{-(\lambda,\alpha_i)}F_{\alpha_i}
    ,\;\forall\lambda\in\Lambda
    \qquad \mbox{\emph{(group action)}}\\
    [E_{\alpha_i},F_{\alpha_j}]
    &=\delta_{i,j}\cdot\frac{K_{\alpha_i}-K_{\alpha_i}^{-1}}
      {v_{\alpha_i}-v_{\alpha_i}^{-1}}
    \qquad \mbox{\emph{(linking)}}
    \end{align*}
    and two sets of \emph{quantum Serre-relations} for any $i\neq j\in I$
    \begin{align*}
    \sum_{r=0}^{1-a_{ij}} (-1)^r
    \begin{bmatrix}1-a_{ij}\\ r\end{bmatrix}_{v_{\alpha_i}}
    E_{\alpha_i}^{1-a_{ij}-r}E_{\alpha_j}E_{\alpha_i}^{r}
    &=0\\
    \sum_{r=0}^{1-a_{ij}} (-1)^r 
    \begin{bmatrix}1-a_{ij}\\ r\end{bmatrix}_{v_{\alpha_i}}
    F_{\alpha_i}^{1-a_{ij}-r}F_{\alpha_j}F_{\alpha_i}^{r}
    &=0
  \end{align*}
  The coproduct, the counit and the antipode are defined on the group algebra
  $\Q(v)[\Lambda]$ as usual  and on the
  additional generators $E_{\alpha_i},F_{\alpha_i}$ as follows:
  \begin{align*}
    \Delta(E_{\alpha_i})
    =E_{\alpha_i}\otimes K_{\alpha_i}+1\otimes E_{\alpha_i} 
    &\qquad \Delta(F_{\alpha_i})
    =F_{\alpha_i}\otimes 1+K_{\alpha_i}^{-1}\otimes F_{\alpha_i}\\
    S(E_{\alpha_i})=-E_{\alpha_i}K_{\alpha_i}^{-1}
    &\qquad S(F_{\alpha_i})=-K_{\alpha_i}F_{\alpha_i}\\
    \epsilon(E_{\alpha_i})=0
    &\qquad \epsilon(F_{\alpha_i})=0
  \end{align*}
 Using Lusztig symmetry operators $T_w$ for any Weyl group element $w$, one can construct \emph{root vectors} $E_\alpha,F_\alpha$ for all $\alpha\in\Phi^+$ such that    
 multiplication in the algebra gives a bijective linear map (PBW-basis): 
 $$\left( \bigotimes_{\alpha\in\Phi^+} \Q(v)[E_\alpha])\right)
 \otimes \Q(v)[\Lambda]
 \otimes \left( \bigotimes_{\alpha\in\Phi^+} \Q(v)[F_\alpha]\right)
 \stackrel{\sim}{\longrightarrow} 
 %U_v^{\Q(v)}^+U_v^{\Q(v)}^0U_v^{\Q(v)}^-=
 U_v^{\Q(v)}$$
 The three subalgebras generated by $E_{\alpha_i},K_{\lambda},F_{\alpha_i}$ respectively are called $U_v^+,U_v^0,U_v^-$.\\
 \item \textbf{Lusztig integral form of divided powers} $U_v^{\Z[v,v^{-1}],\L}(\g)$ (\cite{Lus90} Thm 6.7): An infinite-dimensional Hopf algebra  over the commutative ring  of Laurent 	polynomials $\Z[v,v^{-1}]$,
  generated as a $\Z[v,v^{-1}]$-subalgebra of the rational form by 
  $$E_{\alpha}^{(t)}:=\frac{E_{\alpha}^t}{[t]_{v_\alpha}!},\quad F_{\alpha}^{(t)}:=\frac{F_{\alpha}^t}{[t]_{v_\alpha}!}, \quad   
  K_{\alpha_i},\quad \Ksymb{K_{\alpha}}{0}{t}
  :=\prod_{s=1}^{t}\frac{K_{\alpha}v_{\alpha}^{1-s}
      -K_{\alpha}^{-1}v_{\alpha}^{-1+s}}
      {v_{\alpha}^s-v_{\alpha}^{-s}}$$
  for all $\alpha\in\Phi^+,t\in\N$. Again multiplication in the algebra gives a bijective $\Z[v,v^{-1}]$-linear map (PBW-basis):
  \begin{align*}
 &\left( \bigotimes_{t\in \N,\alpha\in\Phi^+} E_{\alpha}^{(t)}\Z[v,v^{-1}] \right)
 \otimes \left( \bigotimes_{t\leq 1, i}K_{\alpha_i}^t\Z[v,v^{-1}] \right)\\
 &\otimes \left( \bigotimes_{t\in \N}\Ksymb{K_{\alpha}}{0}{t}\Z[v,v^{-1}] \right)
 \otimes \left( \bigotimes_{t\in \N,\alpha\in\Phi^+} F_{\alpha}^{(t)}\Z[v,v^{-1}]\right)
 \stackrel{\sim}{\longrightarrow} 
 U_v^{\Z[v,v^{-1}],\L}   
  \end{align*}
 where we take $K_{\alpha_i}^t,t\leq 1$ instead of all $K_\alpha$ because $K_{\alpha_i}-K_{\alpha_i}^{-1}$ is a multiple of $\Ksymb{K_{\alpha_i}}{0}{1}$.\\
  
  The Lusztig integral form is a Hopf subalgebra of the rational form and has the property that extension of scalars gives an isomorphism of Hopf algebras over $\Q(v)$:
  $$U_v^{\Q(v)}\cong U_v^{\Z[v,v^{-1}],\L}\otimes_{\Z[v,v^{-1}]} \Q(v)$$
  
  \item \textbf{Lusztig quantum group of divided powers} $U_q^{\L}(\g)$ (or restricted specialization):
  An infinite-dimensional Hopf algebra over $\C$ obtained for every choice of an element $q\in\C^\times$, defined by specialization
  $$U_q^{\L}:=U_q^{\Z[v,v^{-1}],\L} \otimes_{\Z[v,v^{-1}]} \C_q$$
  where the indeterminate $v$ acts on $\C_q$ by multiplication with the number $q$. This has the following consequence on the algebra structure: If we have elements in the chosen integral form and a relation between them, including as scalar factors some nonzero $f(v)\in\Z[v,v^{-1}]$, which is clearly invertible in $\Q(v)$, then this may nevertheless specialize to $f(q)=0$. For example, if $q$ a primitive $\ell$-th root of unity we have $[\ell_\alpha]_{q_\alpha}=0$, and correspondingly in $U_q^{\L}(\g)$ we have the relation:
  $$(E_\alpha)^{\ell_{\alpha}}=[\ell_\alpha]!\cdot E_\alpha^{(\ell_{\alpha})}=0,
  \qquad K_{\alpha_i}^{2\ell_{\alpha_i}}=1,\qquad F_\alpha^{\ell_{\alpha}}=0$$
  Multiplication in the algebra gives a bijective linear map (PBW-basis):
  \begin{align*}
   &\left( \bigotimes_{\alpha\in\Phi^+} \C[E_\alpha]/(E_\alpha^{\ell_\alpha}) \otimes \C[E_\alpha^{(\ell_\alpha)}])\right)
 \otimes \C[\Lambda/\Lambda']\\
 &\otimes \left( \bigotimes_{\alpha\in\Phi^+}\C[\Ksymb{K_{\alpha}}{0}{\ell_\alpha}])\right)
 \otimes \left( \bigotimes_{\alpha\in\Phi^+} \C[F_\alpha]/(F_\alpha^{\ell_\alpha})\otimes \C[F_\alpha^{(\ell_\alpha)}])\right)
 \stackrel{\sim}{\longrightarrow} 
 U_q^{\L}
  \end{align*}
   where $\Lambda'$ is the sublattice of the root lattice $\Lambda$ generated by all $K_{\alpha_i}^{2\ell_i}$.\\
   \item \textbf{Small quantum group} $u_q(\g)$: A finite-dimensional Hopf algebra over $\C$ generated as a subalgebra of $U_q^\L(\g)$ by the elements $E_\alpha,K_\lambda,F_\alpha$ for all $\alpha\in\Phi^+$ with 
   $\ell_\alpha >1$, see \cite{Lus90} Sec. 8.2. In the present article we have excluded the small values $\ell=1,2$ and also $\ell=4$ if some $d_\alpha=2$ and $\ell=3,4,6$ if some $d_\alpha=3$. With this assumption in place $\ell_\alpha >1$ holds for all roots $\alpha$ and moreover $u_q(\g)$ is generated by the simple root vectors $E_{\alpha_i},F_{\alpha_i}$ and the $K_\lambda$, otherwise see \cite{Len14}. In particular with this assumption in place $u_q(\g)$ is isomorphic to the small quantum group defined by generators and relations. Multiplication in the algebra gives a bijective linear map (PBW-basis): 
 $$\left( \bigotimes_{\alpha\in\Phi^+} \C[E_\alpha]/(E_\alpha^{\ell_\alpha})\right)
 \otimes \C[\Lambda/\Lambda']
 \otimes \left( \bigotimes_{\alpha\in\Phi^+} \C[F_\alpha]/(F_\alpha^{\ell_\alpha})\right)
 \stackrel{\sim}{\longrightarrow} u_q$$ 
    
  \item We have by \cite{Lus93} Thm. 35.1.9\footnote{I wish to thank David Jordan for pointing out this reference, which covers the general case.} the \emph{Frobenius homomorphism} of algebras:
  $$U_q^{\L}(\g)  \stackrel{\mathrm{Frob}}{\longrightarrow} U(\g^{(\ell)})$$
  by sending the divided powers $E_\alpha^{(\ell_\alpha)},F_\alpha^{(\ell_\alpha)}$ and the expression $\Ksymb{K_{\alpha}}{0}{\ell_\alpha}$ to generators of the following algebra $U(\g^{(\ell)})$: If all $\ell_\alpha$ are equal to $\ell$, or equivalently all $(2d_\alpha,\ell)=1$, then $U(\g^{(\ell)})=U(\g)$ is the universal enveloping of the Lie algebra $\g$. This is the commonly known case treated in \cite{Lus90}. If all $\ell_\alpha$ are equal, or equivalently $4\nmid\ell$ if some $d_\alpha=2$ and $3\nmid\ell$ is some $d_\alpha=3$, then $U(\g^{(\ell)})$ is similar to the universal enveloping of the Lie algebra, with additional signs according to a symmetric braiding $q^{(\ell_{\alpha_i}\alpha_i, \ell_{\alpha_j}\alpha_j)}$. If some $\ell_\alpha$ are not equal, for example in the case where all $2d_\alpha|\ell$, then $U(\g^{(\ell)})$ is similar to the universal enveloping of the Langlands dual Lie algebra, again with possibly additional signs from the symmetric braiding. All of this is also discussed in \cite{Len14}.
  
  We remark, without giving a proof in the most general case (see \cite{A96} for the case $(2d_\alpha,\ell)=1$ and \cite{Len14} for the positive part in the general case) that this should in fact give an exact sequence of Hopf algebras
  $$0\longrightarrow u_q(\g) \longrightarrow U_q^{\L}(\g) \stackrel{\mathrm{Frob}}{\longrightarrow} U(\g^{(\ell)})\longrightarrow 0$$
\end{itemize}
Similarly, the \textbf{Kac-Procesi-DeConcini integral form} $U_q^{\Z[v,v^{-1}],\K}(\g)$ is generated by elements 
 $$E_\alpha^t,\quad F_\alpha^t,\quad  K_\lambda,\quad t\in \N$$
  This (maybe more obvious) choice, which is properly contained in Lusztig's integral form, specializes to a Hopf algebra $U_q^{\K}$ (unrestricted specialization) with a large center properly containing  $E_\alpha^{\ell_\alpha}, F_\alpha^{\ell_\alpha},K_\alpha^{2\ell_\alpha}$. Conjecturally we would expect accordingly for arbitrary $\ell$ an exact sequence:
  %\footnote{This is proven if $\ell$ is relatively prime to the $2d_\alpha$ without $\g^\vee$ appearing, but for in the case of general $\ell$ one would need an analog result to %\cite{Len14} for this integral form and specialization} 
  $$0\longleftarrow u_q(\g) \longleftarrow U_q^{\K}(\g)\longleftarrow \mathcal{O}(G^{(\ell)})\longleftarrow 0$$
  where $\mathcal{O}(G^{(\ell)})$ is a (super-)commutative algebra of functions on the Lie group resp. dual Lie group underlying $\g$ and the choice of $\Lambda$. 
  
  Note that similar to these two specializations one studies two versions of universal enveloping of a Lie algebra over a field in finite characteristic. In the obvious version one gets a large center with more primitive elements, in the non-obvious version one gets additional generators (divided powers) generating a new Lie algebra enveloping, which in this cases however are again truncated.\\   
  
 The Kac-Procesi-DeConcini quantum group will not be used in the main construction of this article, but it is crucial for the application to quantum topology, see Section \ref{sec_ApplicationCategory}. 
 
\section{Main result}

We already claimed that the images of the divided powers $E_\alpha^{(\ell_\alpha)},F_\alpha^{(\ell_\alpha)}$ together with the expression $\Ksymb{K_{\alpha}}{0}{\ell_\alpha}$ map to a basis of the (dual) Lie algebra 
$$U_q^{\L}(\g) \longrightarrow U(\g^{(\ell)})$$
A different more involved preimage of the Chevalley generators of the Cartan part of $U(\g^{(\ell)})$ would be the commutator of the preimages  $[E_\alpha^{(\ell_\alpha)},F_\alpha^{(\ell_\alpha)}]$. Neither preimage is primitive. 
We wish to find such a nicer preimage $H_\alpha$, which already appeared in our proof of \cite{Len14} Thm. 4.1:

\begin{theorem}\label{thm_main}~ 
\begin{enumerate}[a)]
 \item Define the following elements\footnote{Morally the reader should have in mind $H_\alpha=\frac{K_{\alpha}^{2\ell_\alpha}-1}{v_\alpha^{2\ell_\alpha}-1}$, which differs in the specialization by a nonzero scalar.} in the rational form $U_v^{\Q(v)}(\g)$ 
 $$H_{\alpha}':=\frac{K_{\alpha}^{2\ell_\alpha}-1}{\Phi_{\ell_\alpha}(v_{\alpha}^2)}
 %\; \in 
 %U_q^{\Z[v,v^{-1}],\mathcal{L}}
 $$ 
 with $\Phi_{k}(X)$ the $k$-th cyclotomic polynomial. We claim that this element is contained in the integral form $H_\alpha'\in U_q^{\Z[v,v^{-1}],\mathcal{L}}$. More precisely, we claim the following explicit formula  
 %The $H_{\alpha}$ are $\Ksymb{K_{\alpha}}{0}{\ell_\alpha}$ plus additional non-obvious $K_\alpha$-terms made explicit in the proof.
  \begin{align*}
   H_\alpha'=\;v_{\alpha}^{\frac{\ell_{\alpha}(\ell_{\alpha}-1)}{2}}\frac{\prod_{s=1}^{\ell_{\alpha}} v_{\alpha}^s
   -v_{\alpha}^{-s}}{\Phi_{\ell_{\alpha}}(v_{\alpha}^2)}\cdot K_{\alpha}^{\ell_{\alpha}}\Ksymb{K_{\alpha}}{0}{\ell_{\alpha}} 
   \;+\;\frac{(K_{\alpha}^{2\ell_{\alpha}}-1)-(\prod_{s=1}^{\ell_{\alpha}}K_{\alpha}^{2}-v_{\alpha}^{2(s-1)})}{\Phi_{\ell_{\alpha}}(v_{\alpha}^2)}
   \end{align*}
   where the quotient in the first summand is in $\Z[v,v^{-1}]$ and the second summand is a $\Z[v,v^{-1}]$-polynomial in $K_\alpha^2$.
 \item In the specialization $U_q^{\L}(\g)$ with $q$ an $\ell$-th root of unity we define the rescaled element
   $$  H_\alpha
   :=\left(\frac{q_{\alpha}^{2\ell_\alpha}-1}{\Phi_{\ell_{\alpha}}(q_{\alpha}^{2})}\right)^{-1}H_\alpha$$
   where the bracket is understood to be a nonzero value after polynomial division. Then
   \begin{align*}
   H_\alpha
   %&:=\left(\frac{q_{\alpha}^{2\ell_\alpha}-1}{\Phi_{\ell_{\alpha}}(q_{\alpha}^{2})}\right)^{-1}H_\alpha' 
   &=q_\alpha^{\ell_\alpha}(-1)^{\ell_\alpha-1}\ell_\alpha
   \cdot K_{\alpha}^{\ell_{\alpha}}\Ksymb{K_{\alpha}}{0}{\ell_{\alpha}} 
   +\frac{(K_{\alpha}^{2\ell_{\alpha}}-1)-(\prod_{s=1}^{\ell_{\alpha}}K_{\alpha}^{2}-q_{\alpha}^{2(s-1)})}{q_{\alpha}^{2\ell_\alpha}-1}
   \end{align*}
  Moreover we claim that the second summand is the unique $\C$-polynomial of degree $<\ell_\alpha$ in the variable $X=K_\alpha^2$ that takes on each $X=q_\alpha^{2\bar{k}}$ the value $\bar{k}$, where $0\leq \bar{k}<\ell_\alpha$.
 \item The elements $H_{\alpha}'$ in $U_q^{\Z[v,v^{-1}],\mathcal{L}}$ are skew-primitive and have a nice adjoint action
 \begin{align*}
  \Delta(H_{\alpha}')
  &=K_{\alpha}^{2\ell_\alpha}\otimes H_{\alpha}'+H_{\alpha}' \otimes 1\\
  H_{\alpha}'E_{\beta}-v_{\alpha}^{2\ell_\alpha\cdot \frac{2(\alpha,\beta)}{(\alpha,\alpha)}} E_{\beta}H_{\alpha}'
  &=E_{\beta}
  \cdot \frac{v_{\alpha}^{2\ell_\alpha\cdot \frac{2(\alpha,\beta)}{(\alpha,\alpha)}}-1}{\Phi_{\ell_{\alpha}}(v_{\alpha}^{2})}
  %\cdot v_{\alpha}^{2\ell_\alpha (\frac{2(\alpha,\beta)}{(\alpha,\alpha)}-1)} [\frac{2(\alpha,\beta)}{(\alpha,\alpha)}]_{v_{\alpha}^{2\ell_\alpha}} 
  \end{align*}
  \item In the specialization $U_q^{\L}(\g)$ with $q$ an $\ell$-th root of unity, the elements $H_\alpha',H_\alpha$ are primitive and they have the expected adjoint action: 
  \begin{align*}
  \Delta(H_{\alpha})
  &=1\otimes H_{\alpha}+H_{\alpha} \otimes 1\\
  H_{\alpha}E_{\beta}-E_{\beta}H_\alpha
  &=\frac{2(\alpha,\beta)}{(\alpha,\alpha)}E_{\beta}
  %\cdot \frac{v_{\alpha}^{2\ell_\alpha\cdot \frac{2(\alpha,\beta)}{(\alpha,\alpha)}}-1}{\Phi_{\ell_{\alpha}}(v_{\alpha}^{2\ell_\alpha})}
  %\cdot v_{\alpha}^{2\ell_\alpha (\frac{2(\alpha,\beta)}{(\alpha,\alpha)}-1)} [\frac{2(\alpha,\beta)}{(\alpha,\alpha)}]_{v_{\alpha}^{2\ell_\alpha}} 
  \end{align*}
  %where after polynomial division the evaluation $\frac{q_{\alpha}^{2\ell_\alpha}-1}{\Phi_{\ell_{\alpha}}(q_{\alpha}^{2})}$ is a well-defined nonzero complex number. It can be %removed by rescaling $H_\alpha$.
  \item We hence find the \emph{unrolled small quantum group} inside $U^{\L}_q(\g)$: 
   $$u_q(\g) {U^{\L}_q(\g)}^0
   \cong u_q(\g)\rtimes U(\h)$$
   with the image of $\h$ given by the $H_\alpha$.\\
\end{enumerate}
\end{theorem}
Before we prove the theorem, we discuss the result: 
\begin{example}
    For $\g=\sl_2$ and $\ell=4$ we have $d_\alpha=1,\ell_\alpha=2,q_\alpha^{\ell_\alpha}=-1$ for the single root $\alpha$, which we from now on omit from the notation. Then the results of the theorem read:
    \begin{align*}
      H'
      &:=\frac{K^4-1}{v^2+1}
      =v\frac{(v-v^{-1})(v^2-v^{-2})}{v^2+1}K^2\Ksymb{K}{0}{2}
      +(K^2-1) \\
      H&:=\left(\frac{q^4-1}{q^2+1}\right)^{-1}H'
      %=q\frac{(q-q^{-1})(q^2-q^{-2})}{q^4-1}K^2\Ksymb{K}{0}{2}
      %+\frac{K^2-1}{q^2-1} 
      =2K^2\Ksymb{K}{0}{2}
      +\frac{K^2-1}{-2}
    \end{align*}
    and as asserted the second summand $\frac{K^2-1}{-2}$ is equal to $0$ or $1$, respectively, if we set $K^2=q^{2k}$ with $k\equiv 0$ or $k\equiv 1$ modulo $2$, respectively.
\end{example}
\begin{corollary}
    Theorem \ref{thm_main} e) shows that the unrolled quantum group appearing in \cite{GPT09,CGP15,GP16} modulo the relation $K_{\alpha_i}^{2\ell_i}=1$ appears as Hopf subalgebra of $U_q^{\L}(\g)$, see Section~\ref{sec_ApplicationCategory}.  
\end{corollary}
\begin{corollary}\label{cor_weightshift}
    Theorem \ref{thm_main} b) has the following direct application: Let $W$ be a module of $U_q^\L(\g)$ and $w_\lambda$ a highest-weight vector for $u_q(\g)$ in the sense that for every simple root $\alpha$
    $$E_{\alpha}w_\lambda=0,\qquad
    K_\alpha w_\lambda=q_\alpha^{(\alpha,\lambda)}=q^{(\alpha^\vee,\lambda)}w_\lambda, \qquad 
    H_\alpha w_\lambda =(\alpha^\vee,\lambda)w_\lambda$$
    In $U_q^\L(\g)$ we have from \cite{Lus90} 6.5(a2) the following commutator relation for simple root vectors
  $$\left[E_{\alpha}^{(\ell_{\alpha})},F_{\alpha}^{(\ell_{\alpha})}\right]
  =\Ksymb{K_{\alpha}}{0}{\ell_{\alpha}}+\sum_{t=1}^{\ell_{\alpha}-1} F_{\alpha}^{(\ell_{\alpha}-t)} \Ksymb{K_{\alpha}}{2t-2\ell_{\alpha}}{t}  E_{\alpha}^{(\ell_{\alpha}-t)}$$
  so the commutator acts on $w_\lambda$ in the same way as the divided $K$-power. \\
  
  Theorem \ref{thm_main} b) calculates this action: Define $k=(\alpha^\vee,\lambda)$ and write $k=\bar{k}+k'\ell_\alpha$ with $0\leq \bar{k} < \ell_{\alpha}$, then the second summand in b) acts on $w_\lambda$ by multiplication with $\bar{k}$, and since $H_\alpha$  acts by multiplication with $k$, the formula in b) implies further 
  $$\left[E_{\alpha}^{(\ell_{\alpha})},F_{\alpha}^{(\ell_{\alpha})}\right]w_\lambda=\Ksymb{K_{\alpha}}{0}{\ell_\alpha}w_\lambda= q_\alpha^{\ell_\alpha}(-1)^{\ell_\alpha-1}k'\cdot w_\lambda  $$
  More Lie-theoretically spoken, $H_\alpha$ and the commutator $[E_{\alpha}^{(\ell_{\alpha})},F_{\alpha}^{(\ell_{\alpha})}]$ act on such vectors annihilated by $E_\alpha$ by weights $\lambda$ and $\ell_\alpha\lambda'$, which essientially differ by the unique representative $\bar{\lambda}$ in a fundamental alcove. This shift also appears on the CFT side, see Section~\ref{sec_ApplicationCFT}.
\end{corollary}~

\begin{proof}[Proof of Theorem \ref{thm_main}]~
 \begin{enumerate}[a)]
  %%%%%%%%%%%%%%%%%%%%%%
  \item 
The following similar element from \cite{Lus90} is by definition clearly contained in $U_q^{\Z[v,v^{-1}],\mathcal{L}}$  
  \begin{align*}
  &K_{\alpha}^{\ell_{\alpha}}\Ksymb{K_{\alpha}}{0}{\ell_{\alpha}}
  =K_{\alpha}^{\ell_{\alpha}}
      \prod_{s=1}^{\ell_{\alpha}}\frac{K_{\alpha}v_{\alpha}^{1-s}
      -K_{\alpha}^{-1}v_{\alpha }^{s-1}}
      {v_{\alpha}^s-v_{\alpha}^{-s}}
  =v_{\alpha}^{-\frac{\ell_{\alpha}(\ell_{\alpha}-1)}{2}}\cdot
      \frac{\prod_{s=1}^{\ell_{\alpha}} K_{\alpha}^{2}-v_{\alpha}^{2(s-1)}}
      {\prod_{s=1}^{\ell_{\alpha}} v_{\alpha}^s-v_{\alpha}^{-s}}, 
   \end{align*}
   Subtracting a suitable $\Z[v,v^{-1}]$-rescaling of this element from our element $H_{\alpha}$ gives
   \begin{align*}
   \frac{K_{\alpha}^{2\ell_{\alpha}}-1}{\Phi_{\ell_{\alpha}}(v_{\alpha}^2)} 
   \;-\;v_{\alpha}^{\frac{\ell_{\alpha}(\ell_{\alpha}-1)}{2}}\frac{\prod_{s=1}^{\ell_{\alpha}} v_{\alpha}^s
   -v_{\alpha}^{-s}}{\Phi_{\ell_{\alpha}}(v_{\alpha}^2)}\cdot K_{\alpha}^{\ell_{\alpha}}\Ksymb{K_{\alpha}}{0}{\ell_{\alpha}} 
   &=\frac{(K_{\alpha}^{2\ell_{\alpha}}-1)-(\prod_{s=1}^{\ell_{\alpha}}K_{\alpha}^{2}-v_{\alpha}^{2(s-1)})}{\Phi_{\ell_{\alpha}}(v_{\alpha}^2)}
   \end{align*}
   This is a polynomial in $K_{\alpha}^2,v_{\alpha}^2$, because the numerator has a zero at $v_{\alpha}^2=q_{\alpha}^2$, so 
   if we write the numerator as $\sum_n c_n(v_{\alpha}^2) K_{\alpha}^{2n}$ then each polynomial $c_n(v_{\alpha}^2)$ is divisible by $\Phi_{\ell_{\alpha}}(v_{\alpha}^2)$.\\

  \item Before specialization we rewrite in the rational form
  \begin{align*}
   H_\alpha
   &:=\left(\frac{v_{\alpha}^{2\ell_\alpha}-1}{\Phi_{\ell_{\alpha}}(v_{\alpha}^{2})}\right)^{-1}H_\alpha' \\
   &=v_{\alpha}^{\frac{\ell_{\alpha}(\ell_{\alpha}-1)}{2}}\frac{\prod_{s=1}^{\ell_{\alpha}} v_{\alpha}^s
   -v_{\alpha}^{-s}}{v_{\alpha}^{2\ell_\alpha}-1}\cdot K_{\alpha}^{\ell_{\alpha}}\Ksymb{K_{\alpha}}{0}{\ell_{\alpha}} 
   +\frac{(K_{\alpha}^{2\ell_{\alpha}}-1)-(\prod_{s=1}^{\ell_{\alpha}}K_{\alpha}^{2}-v_{\alpha}^{2(s-1)})}{v_{\alpha}^{2\ell_\alpha}-1}\\
   &=v_\alpha^{\ell_\alpha}(-1)^{\ell_\alpha-1}
   \prod_{s=1}^{\ell_{\alpha}-1} (1-v_{\alpha}^{2s})
   \cdot K_{\alpha}^{\ell_{\alpha}}\Ksymb{K_{\alpha}}{0}{\ell_{\alpha}} 
   +\frac{(K_{\alpha}^{2\ell_{\alpha}}-1)-(\prod_{s=1}^{\ell_{\alpha}}K_{\alpha}^{2}-v_{\alpha}^{2(s-1)})}{v_{\alpha}^{2\ell_\alpha}-1}\\
   \end{align*}
   The value of the first summand at $v=q$ is up to a factor the polynomial $(X^{\ell_\alpha}-1)/(X-1)$ at $X=1$, which is  equal $\ell_\alpha$, and hence the first summand has value 
   $$=q_\alpha^{\ell_\alpha}(-1)^{\ell_\alpha-1}\ell_\alpha
   \cdot K_{\alpha}^{\ell_{\alpha}}\Ksymb{K_{\alpha}}{0}{\ell_{\alpha}}$$
   For the second assertion we consider the second summand before specialization as a $\Q(v)$-polynomial in the variable $X=K_\alpha^2$, which is clearly of degree $<\ell_\alpha$ because the highest power cancels. For any $0\leq \bar{k}<\ell_\alpha$ we observe that for $X=v_\alpha^{\bar{k}}$ the product has a zero at $s=\bar{k}$, thus
   $$\frac{(X^{\ell_{\alpha}}-1)-(\prod_{s=1}^{\ell_{\alpha}}X^{2}-v_{\alpha}^{2(s-1)})}{v_{\alpha}^{2\ell_\alpha}-1}
   =\frac{v_\alpha^{2\bar{k}\ell_{\alpha}}-1}{v_{\alpha}^{2\ell_\alpha}-1}$$
   It is clear that this specializes to $\bar{k}$. For the uniquess assertion, we denote that any two complex polynomials with prescribed values at all $\ell_\alpha$-th roots of unity $X=(q_\alpha^2)^{\bar{k}}$ have as difference a multiple of $X^{\ell_\alpha}-1$. Hence there is at most one such polynomial of degree $<\ell_\alpha$.
 
  %%%%%%%%%%%%%%%%%%%%%%
  \item The desired relations hold in the rational form $U_v^{\Q(v)}$:
\begin{align*}
  \Delta(H_{\alpha}')
  &=\frac{\Delta(K_\alpha^{2\ell_i}-1)}{\Phi_{\ell_\alpha}(v_\alpha^2)}\\
  &=\frac{K_{\alpha}^{2\ell_\alpha}\otimes (K_\alpha^{2\ell_i}-1)+(K_\alpha^{2\ell_i}-1) \otimes 1}{\Phi_{\ell_\alpha}(v_\alpha^2)}\\
  &=K_{\alpha}^{2\ell_\alpha}\otimes H_{\alpha}+H_{\alpha} \otimes 1\\
   H_{\alpha}'E_{\beta}-v_{\alpha}^{2\ell_\alpha\cdot \frac{2(\alpha,\beta)}{(\alpha,\alpha)}} E_{\beta}H_{\alpha'}
  &=E_{\beta}
  \left(\frac{K_{\alpha}^{2\ell_\alpha}\left(v_{\alpha}^{\frac{2(\alpha,\beta)}{(\alpha,\alpha)}}\right)^{2\ell_\alpha}-1}{\Phi_{\ell_{\alpha}}(v_{\alpha}^{2})}
  -v_{\alpha}^{2\ell_\alpha\cdot \frac{2(\alpha,\beta)}{(\alpha,\alpha)}}\frac{K_{\alpha}^{2\ell_\alpha}-1}{\Phi_{\ell_{\alpha}}(v_{\alpha}^{2\ell_\alpha})}\right)\\
  &=E_{\beta}
  \cdot \frac{v_{\alpha}^{2\ell_\alpha\cdot \frac{2(\alpha,\beta)}{(\alpha,\alpha)}}-1}{\Phi_{\ell_{\alpha}}(v_{\alpha}^{2})}
  %\cdot v_{\alpha}^{2\ell_\alpha (\frac{2(\alpha,\beta)}{(\alpha,\alpha)}-1)} [\frac{2(\alpha,\beta)}{(\alpha,\alpha)}]_{v_{\alpha}^{2\ell_\alpha}} 
  \end{align*}
  \item The last $q$-factor in c) may again be written as $q$-number and specializes for $q_\alpha^{2\ell_\alpha}=1$ to 
$$\frac{2(\alpha,\beta)}{(\alpha,\alpha)}\cdot 
\frac{v_\alpha^{2\ell_\alpha}-1}{\Phi_{\ell_{\alpha}}(v_{\alpha}^{2})}
$$
  On the other hand the $K_\alpha^{2\ell_\alpha}$ in the coproduct specializes to $1$.
  %%%%%%%%%%%%%%%%%%%%%
  \item The unrolled quantum group is defined in \cite{GP16} as a smash-product $u_q(\g)\rtimes U(\h)$ with $\h$ having basis $h_{\alpha_i}$, and the action of $U(\h)$ on $u_q(\g)$ is given by rescaling:
  $$h_{\alpha_i}.E_{\alpha_j}=(\alpha_i^\vee,\alpha_j)E_{\alpha_j},\qquad
  h_{\alpha_i}.F_{\alpha_j}=-(\alpha_i^\vee,\alpha_j)F_{\alpha_j},\qquad
  h_{\alpha_i}.K_{\lambda}=0$$
  which gives the following algebra relations in the smash product
  $$[h_{\alpha_i},E_{\alpha_j}]=(\alpha_i^\vee,\alpha_j)E_{\alpha_j},\qquad
  [h_{\alpha_i},F_{\alpha_j}]=-(\alpha_i^\vee,\alpha_j)F_{\alpha_j},\qquad
  [h_{\alpha_i},K_{\lambda}]=0$$
  We first prove that there is a Hopf algebra homomorphism 
  $$f:\;u_q(\g)\rtimes U(\h){\longrightarrow} U_q^{\L}(\g)$$
  As linear map between the underlying vector spaces $f:u_q(\g)\otimes U(\h)\to U_q^\L(\g)$   
  we define $f(a\otimes h_{\alpha_i})\mapsto a\cdot H_{\alpha_i}$ with $H_{\alpha_i}$ defined in b). This is an Hopf algebra map on the Hopf subalgebras $u_q(\g),U(\h)$, because d) shows that the $H_{\alpha_i}$ are primitive. Moreover d) shows that the commutation relations in the smash product match the commutation relations in $U^\L_q(\g)$.\\
  
  To show that the Hopf algebra map $f$ is bijective we use the PBW-basis in $u_q(\g)$ and $U_q^\L(\g)$. Hence in the smash-product every element is a unique linear combination of sorted monomials in $E_{\alpha},F_{\alpha},\alpha\in\Phi^+$ and $K_\lambda$ and $h_{\alpha_i}$, and in $U_q^\L(\g)$ every element is a unique linear combination of sorted monomials in $E_{\alpha},F_{\alpha},\alpha\in\Phi^+$ and $K_\lambda$ and divided powers $E_{\alpha}^{(\ell_{\alpha})},F_{\alpha}^{(\ell_{\alpha})}, \Ksymb{K_{\alpha}}{0}{\ell_\alpha}$. By definition, the image of a monomial in $h_{\alpha_i}$ under $f$ is a nonzero multiple of the respective monomial in $\Ksymb{K_{\alpha_i}}{0}{\ell_{\alpha_i}}$ plus terms involving less divided $K$-powers. Hence, if we enumerate the PBW-basis of the vector space  $U_q^\L(\g)$ by increasing number of divided $K$-power factors, then $f$ is an upper triangular matrix, with nonzero diagonal entries assigning to a PBW monomial 
  in $E_{\alpha},F_{\alpha},K_\lambda,h_{\alpha_i}$ the respective PBW monomial in $E_{\alpha},F_{\alpha},K_\lambda,\Ksymb{K_{\alpha_i}}{0}{\ell_{\alpha_i}}$. In particular $f$ is a bijection and the image is equal to $u_q(\g)U_q^\L(\g)^0$ as asserted.
 \end{enumerate}
\end{proof}

\section{Unrolled quantum groups and categories}\label{sec_ApplicationCategory}

In this section we will show the applications of our results to the constructions in \cite{GPT09,CGP15,GP16} by which it has been inspired and for which it has been intended:\\

In these articles, the authors construct quantum invariants for manifolds and knots from nonsemisimple categories by ``modified traces'', which solves the problem that invariants obtained from nonsemisimple tensor categories tend to vanish. The first example in mind might be representations of a small quantum group $u_q$, but for an even root of unity this does not have a braiding.\\

To get a modular tensor category, they instead consider the unrestricted quantum group $U_q^{\K}$, and then impose relations $E^r=F^r=0$ (this is dropped in \cite{GP16}) while $K$ remains of infinite order; let us call this quotient $\bar{U}_q^{\K}$. Then they unroll this quantum group to $\bar{U}_q^{\K}\rtimes U(\h)$ and consider the category of modules $\mathcal{C}$ where $K$ acts as $q^h,h\in\h$.\\ 

The author believes that there are two disjoint mechanisms at work here: 
\begin{itemize}
 \item The unrolling seems to be intended to get a modular tensor category, because it removes ambiguities in the weight $\lambda$ associated to a $K_\alpha$-eigenvalue $q_\alpha^{(\alpha,\lambda)}$. This should work and apply already in the case of the unrolled \emph{small} quantum group treated in this article.\\
 
 Besides unrolling there are several attempts to construct modular tensor categories for small quantum groups at even roots of unity. The author would assume that they all can in principle be used for the same purpose:
 \begin{itemize}
 \item Accept that the braiding is only up to an outer automorphism as in \cite{Tan92,Res95}.
  \item Consider larger lattices $\Lambda$ and then accept that the monodromy matrix is only factorizable on a sub-Hopf algebra as in \cite{RT91}. This seems like a small-scale unrolling, just enough to remove the $\Z_2$-obstruction
  \item Accept the the quadratic form $q^{(\alpha,\alpha)}$ on the finite group ring does not come from a bimultiplicative form, and construct accordingly a quasi-Hopf algebra version of $u_q(\g)$ at even $\ell$, which again admits a braiding and has a modular tensor category of representations as done for $\sl_2$ explicitly in \cite{GR15}. More calculationally, this amounts to choose representatives $\lambda,\mu,\ldots$ for all cosets, in order to be able to get a well-defined braiding $q^{(\lambda,\mu)}$ and accept that this comes at the expense of introducing a nontrivial associator. The same quandratic form and corresponding $3$-cocycle are also clearly visible on the CFT side for the lattice VOA.
  
  We have given a systematic construction of the quasi-Hopf algebra with modular tensor category in \cite{GLO18}. Essentially we construct quantum groups, where already the Cartan part is a quasi-Hopf algebra determined by the quadratic form on the finite abelian group. The resulting quasi-Hopf algebra can also be viewed as modularization of the extended Hopf algebra in the previous bullet.
 \end{itemize}~\\
 
 \item Taking now the unrestricted quantum group $U_q^{\K}$ gives a much larger category. which is fibred over $G^{(\ell)}$ according to the action of $\mathcal{O}(G^{(\ell)})$, e.g. according to the action of $K_\alpha^{2\ell_\alpha}$. In the case $(\ell,2d_\alpha)=1$ the fibres are known to be semisimple categories on elements in regular conjugacy classes of $G^{(\ell)}=G$, for arbitrary $\ell$ we would expect something similar. On the other hand there are singular points, most importabtly the fiber $\mathcal{C}_0$ with all $K_\alpha^{2\ell_\alpha}=1$, where we get the representations of the small quantum group. Studying quantum invariants attached to this situations is fascinating, in particular studying the invariants at the singular points using complex analysis on the variety $G$.\\
 
%  The author's construction can also be applied to the unrolling of this situation. The authors suggestion would moreover be that this unrolled is again a subalgebra of a larger Hopf algebra, which is a curious hybrid of 
%  $U_q^{\K}$ and $U_q^{\mathcal{L}}$. Both points are made explicit in what follows.
 \end{itemize}
 
 The main result of this article realizes the unrolled small quantum group $u_q\rtimes U(\h)$ inside Lusztig's quantum group. Moreover, it makes the condition $K=e^h$ very natural and thus explains the category $\mathcal{C}_0$ from this point of view. \\

 The authors suggestion would moreover be that this unrolled is again a subalgebra of a larger Hopf algebra, which is a curious hybrid of 
 $U_q^{\K}$ and $U_q^{\mathcal{L}}$:
 
 \begin{question}
   Does there exist for any $q\in\C^\times$ an infinite-dimensional complex Hopf algebra $U_q^{\K\L}(\g)$ with the following properties:
 \begin{enumerate}[a)]
  \item There are Hopf algebra maps
 $$U_q^{\K}\hookrightarrow U_q^{\L\K}\qquad U_q^{\K\L}\twoheadrightarrow U_q^{\L}$$
 such that their composition sends $E_\alpha,F_\alpha,K_\lambda$ to themselves. %Probably $U_q^{\K\L}$ is in some sense universal with this property. 
  \item There are elements such that multiplication induces a bijection of vector spaces:
  $$\left( \bigotimes_{\alpha\in\Phi^+} \C[E_\alpha]\otimes \C[E_\alpha^{(\ell_\alpha)}])\right)
 \otimes \C[\Lambda]\otimes \left( \bigotimes_{\alpha\in\Phi^+}\C[H_\alpha])\right)
 \otimes \left( \bigotimes_{\alpha\in\Phi^+} \C[F_\alpha]\otimes \C[F_\alpha^{(\ell_\alpha)}])\right)
 \stackrel{\sim}{\longrightarrow} 
 U_q^{\K\L}(\g)$$  
 \end{enumerate}
 \end{question}
 
We like to remark that \cite{GG17} Sec. 5.4 have indeed constructed a ``mixed quantum group'', where the negative part resembles $U^\K$ and the positive part resembles $U^\L$.

\newpage
\section{Diagonal Nichols algebras}\label{sec_ApplicationNichols}

We argue, that the same construction should be done for any Nichols algebra of diagonal type: \\

Let $(V,c)$ be a vector space with a braiding $c:V\otimes V\to V\otimes V$, then the Nichols algebra $\B(V)$ is a certain quotient of the tensor algebra of $V$. In this section we just consider diagonal braidings  $c(x_i\otimes x_j)=q_{ij}(x_j\otimes x_i)$ for some basis $x_1,\ldots, x_n\in V$ and arbitrary numbers $q_{ij}\in\C^\times$. In this case we refer the reader to \cite{Heck09}, for the more general case see \cite{AHS10}. \\

\newcommand{\sym}{\mathrm{Sym}}
\begin{definition}
$\B(V)$ is the tensor algebra $\bigoplus_{n\geq 0}V^{\otimes n}$ divided in each degree $n$ by the kernel of the following \emph{quantum symmetrizer map}
$$\sym_{q,n}:\; V^{\otimes n }\to V^{\otimes n} \qquad\qquad  
\sym_{q,n}: =\sum _{\tau\in \mathbb{S}_n} \rho_n(s(\tau))$$
Here, $\tau$ runs over all $n$-permutations, $s(\tau)$ is (any) preimage of $\tau$ of shortest length in the group $\mathbb{B}_n$ of braids on $n$ strands (this set-theoretic section of $\mathbb{B}_n\twoheadrightarrow \mathbb{S}_n$ is called Matsumoto section) and $\rho_n$ is the representation of the braid group on $V^{\otimes n}$ using the given braiding $c_{V,V}:V\otimes V\to V\otimes V$.
\end{definition}

But the Nichols algebra in fact carries much more structure, and most deeper structural results make use of this structure: Assume $V$ is an object in some braided tensor category $\mathcal{C}$ with braiding $c=c_{V,V}$. Then the tensor algebra becomes a Hopf algebra in $\mathcal{C}$, when we define the coproduct by taking $V$ to be primitive elements. The Nichols algebra $\B(V)$ is a graded Hopf algebra inside the braided tensor category $\mathcal{C}$, namely the unique quotient of the tensor algebra with one of several equivalent universal properties: For example, $\B(V)$ is the unique graded Hopf algebra in $\mathcal{C}$ where $V$ is precisely the set of all primitive elements in $\B(V)$. Because of this universal property, Nichols algebras are a crucial ingredient in the classification of certain classes of Hopf algebras in \cite{AS10}.

\begin{example} Let $\g$ be a complex semisimple finite-dimensional Lie algebra of rank $n=\dim(\h)$. Let $q\in\C^\times$ and define the braided vector space
$$V:=\bigoplus_{i=1}^n \C E_{\alpha_i}, \qquad q_{ij}:=q^{(\alpha_i,\alpha_j)}$$
Then $\B(V)$ is isomorphic to $U_q^+(\g)$ if $q$ is not a root of unity \cite{Lus93}, although Nichols algebras are not defined explicitly, and $\B(V)$ is isomorphic to $u_q(\g)^+$ if $q$ is a root of unity \cite{AS00}.
\end{example}
Other examples include the super-Lie algebras with some $q_{ii}=-1$ (fermionic) and some other cases with $q_{ii}$ of small order. A complete classification of finite-dimensional Nichols algebras of diagonal type and a striking structure theory by arithmetic root systems and Weyl groupoids has been given by Heckenberger \cite{Heck09}. The theory of Weyl groupoids has been generalized to a large class of Nichols algebras of non-diagonal type in \cite{AHS10}. 

\begin{example}[$\sl(2|1)$] Let $q\in\C^\times$ and choose 
  $$V:=\C E_1\oplus \C E_2, \qquad q_{ij}:=\begin{pmatrix}
					  -1 & q^{-1}\\
					  q^{-1} & -1
                                       \end{pmatrix}$$
  This datum gives rise to the Nichols algebra $\B(V)=U_q(\sl(2|1))^+$ respectively $\B(V)=u_q(\sl(2|1))^+$, and it makes very transparent how a Weyl group is generalized to a Weyl groupoid: Already for super-Lie algebras there are different non-isomorphic Borel parts, and correspondingly different types of Weyl chambers. Correspondingly, there are some reflections (called odd reflections for super Lie algebras) that changes our initial braiding matrix to a different type of Weyl chamber:  
  $$V=\C F_1\oplus \C E_{12}, \qquad q_{ij}=\begin{pmatrix}
					  -1 & q^{-1}\\
					  q^{-1} & q
                                       \end{pmatrix}$$
   In this example, these two different Weyl chambers still have the same Cartan matrix, but already the super Lie algebra $D(2,1;\alpha)$ demonstrates that this is not always the case. All Weyl groupoids and root systems were classified in \cite{CH09}. 
\end{example}
Let $V$ be as usual realized as a Yetter-Drinfeld module over an abelian group $G$, say $G=\Z^n$, so $q_{ij}$ is given by a bicharacter $\chi(g_i,g_j)$. The quantum double construction can be used in the same way for an arbitrary Nichols algebra of diagonal type to define an analog of the quantum group $U(\chi)$, see \cite{Heck10,AY13}.\\

Now let us consider $U(\chi(v))$ with braiding matrix $q(v)_{ij}=\chi(v)(g_i,g_j)$ depending on a free parameter $v\in\C^\times$ as above, i.e. over the field $\Q(v)$. We can of course proceed precisely as in \cite{Lus90}: 
Use the generalizations of Lusztig's symmetries to construct root vectors ${E_\alpha,F_\alpha\in U(V,q_{ij})}$ for all roots $\alpha$ as in \cite{Heck10}. Then we can define the Hopf subalgebra and $\Z[v,v^{-1}]$-submodule $A$ generated by all PBW monomials in the root vectors resp. by all PBW monomial in divided powers of the root generators. 
 
 Then for a specific value $q\in \C^\times$, e.g. a root of unity, specialization $\otimes_{\Z[v,v^{-1}]} \C_q$ given then two complex Hopf algebras, which could be seen as generalized versions of the Kac-Procesi-DeConcini quantum group $U^\K(\chi)$ and the Lusztig quantum group of divided powers $U^\L(\chi)$.  

\begin{question}
 Does one obtain in this way an integral form 
 $A\otimes_{\Z[v,v^{-1}]} \Q(v)\cong U(\chi)$, so that the PBW-monomials are a vector space basis of the specialization? What can be said about the associated two versions of quantum groups? (for example in relation to the distinguished pre- resp. post-Nichols algebra in \cite{A15} resp. \cite{AAB15}) Moreover, we would also again like to study the hybrid quantum group $U^{\K\L}(\chi)$.
 \end{question}
 Note that the divided powers that appear are only over those root vectors $E_\alpha$ with no a-priori truncation relations in place, e.g. in $\sl(2|1)$ only for $E_{12}$ but not for the fermionic generators with $E_1^2=E_2^2=0$. These are the Cartan-like roots, see e.g. \cite{A15}. The set of Cartan like roots (scaled by their truncation power $\ell_\alpha$) form an ordinary root system, typically of smaller rank than $V$. If we call the respective ordinary Lie algebra $\g^{(\ell)}$ and the corresponding Lie group $G^{(\ell)}$, for example $\g^{(\ell)}=\sl_2$ for $\g=\sl(2|1)$ above, then we again expect exact sequences of Hopf algebras 
$$0\longrightarrow U_q(\chi(q))\longrightarrow U^\L(\chi) \longrightarrow U(\g^{(\ell)}) \longrightarrow 0$$ 
$$0\longleftarrow U_q(\chi(q))\longleftarrow U^\K(\chi) \longleftarrow \mathcal{O}(G^{(\ell)}) \longleftarrow 0$$ 
where $u_q(\chi)=U_q(\chi(q))$ is the (finite-dimensional) Hopf algebra associated to $\chi(q)$.\\
 
 \begin{question}
 Can one realize the unrolled quantum group of $U(\chi)$ in \cite{AS17} inside this $U^\L(\chi)$?
 
 On the other hand, can one construct a modular tensor category of representations of some quasi-Hopf algebra variant of $U(\chi)$? As for the quantum groups at even roots of unity, this probably requires to realize $V$ over a quasi-Hopf algebra associated to a finite abelian group with a quadratic form as in \cite{GLO18}.  
 \end{question}
 \begin{question}
  Are there interesting new invariants of $3$-manifolds attached to this setting? Are there corresponding vertex algebras with this modular tensor category as their category of representation? 
 \end{question}

\section{Conformal field theory}\label{sec_ApplicationCFT}
 
We also wish to point out the connection of our results on unrolling via Lusztig's quantum group to logarithmic conformal field theory \cite{FGSTsl2,NT11, FT10,Len17}:\\

Fix $\g$ and $\ell$ such that all $2d_\alpha|\ell$ (i.e. $\ell_\alpha=\ell/2d_\alpha$) and define $p=\ell/2$. Let $\Lambda=\sqrt{p}\Lambda_{R^\vee}$ a rescaling of the root lattice of $\g^\vee$. Then it has been conjectured that: 
\begin{itemize}
 \item There exists an action of parts of Lusztig quantum group $U^\L_q(\g)$ on the lattice vertex operator algebra $\mathcal{V}_\Lambda$ associated to $\Lambda$, which physically describes a compactified free boson.
 \item The kernel $\mathcal{W}\subset \mathcal{V}_\Lambda$ of the action of the subalgebra $u_q(\g)^+$ has as category of representations a non-semisimple modular tensor category, which is close to the category of $u_q(\g)$-modules. More precisely it is expected to be equivalent as braided tensor category to the category of representations of a quasi-Hopf algebra, possibly the one we recently construct in \cite{GLO18}. 
\end{itemize}
The program has been carried out for $\g=\sl_2$ in \cite{FGSTsl2,NT11}, in the case $p=2$ the quasi-Hopf algebra is obtained in \cite{GR15}. The quantum group relations for the action of $u_q(\g)^+$ have been proven in general by the author in \cite{Len17}.\\

The action of $F_\alpha$ are given by \emph{short screening charge operators} $\res(\Y(\exp{-\alpha/\sqrt{p}}))$, the action of $E_\alpha^{(\ell_\alpha)}$
 by \emph{long screening charge operators} $\res(\Y(\exp{+\alpha^\vee\sqrt{p}}))$, the action of $H_{\alpha}$ by a \emph{scalar charge operator} $\res(\Y(\partial\phi_{\alpha^\vee\sqrt{p}}))$ and the action of $K_\alpha$ by it's rescaled exponential $e^{\pi i\;\res(\Y(\partial\phi_{\alpha/\sqrt{p}}))}$. 
 The evaluation of the scalar charge operator $H_\alpha$ on a element in some module $v_{\lambda}$ with degree $\lambda/\sqrt{p}\in \Lambda^*$ is  
 \begin{align*}
H_{\alpha} v_{\lambda}&=\left(\alpha^\vee\sqrt{p},\lambda/\sqrt{p}\right)=\frac{2(\alpha,\lambda)}{(\alpha,\alpha)}v_\lambda\\
K_\alpha v_\lambda&=e^{\pi i(\alpha/\sqrt{p},\lambda/\sqrt{p})}v_\lambda=e^{\frac{2\pi i}{2p}(\alpha,\lambda)}v_\lambda  
 \end{align*}
 We see that this matches (up to a rescaling due to the dual) the condition $K=q^H$ on the category of the unrolled small quantum group. In this setting we can also recover our formula for $H_\alpha$ from our Main Theorem \ref{thm_main}:  
 \begin{align*}
  H_\alpha v_\lambda
 &= \lim_{v\to q}
 \frac{K_\alpha^{2\ell_\alpha}-1}{v_\alpha^{2\ell_\alpha}-1}v_\lambda 
 =\lim_{v\to q}
 \frac{v^{\frac{2(\alpha,\lambda)}{(\alpha,\alpha)}\ell}-1}{v^{\ell}-1}v_\lambda
 =\frac{2(\alpha,\lambda)}{(\alpha,\alpha)} v_\lambda
 \end{align*}

From a lattice vertex algebra perspective it is natural that a braiding in $\mathcal{V}_\Lambda$ comes not easy, because any module $\mathcal{V}_{[\lambda/\sqrt{p}]}$ contains an entire coset $[\lambda/\sqrt{p}]\in\Lambda^*/\Lambda$ and usually the braiding $q^{(\lambda,\mu)}$ does not factorize over this quotient, so one has to choose representatives, which is general causes a $3$-cocycle to appear as associator. Note that a coset $[\lambda/\sqrt{p}]$ is precisely the set of degrees with the same action of all $K_\alpha$. The effect of ``unrolling'' is to separate the different elements in the coset and hence to be able to define the braiding without  ambiguities. This is the same effect we want for unrolling quantum group representations, and it is to expect the two versions are closely related.   

\begin{observation}
 We observe that the weight shift in our main formula's Corollary \ref{cor_weightshift} by the unique representative $\bar{\lambda}$ in the fundamental alcove is the same shift that appears in \cite{FT10} formula (4.25) on the level of CFTs. Here we are in the divisible case $2d_\alpha | \ell$, so the prefactor $q_\alpha^{\ell_\alpha}(-1)^{\ell_\alpha-1}=(-1)^{\ell_\alpha}$.
 This is another indication, that we have an action of the Lusztig divided power algebra in this CFT setting, not simply an action of $u_q(\g)$ and $U(\g^\vee)$.
\end{observation}

 \begin{acknowledgementX}
  The author thanks C. Schweigert for interesting discussions, input and support. 
  The author receives additional support by the DFG Graduiertenkolleg RTG 1670 at the University of Hamburg. 
 \end{acknowledgementX}

\end{document}